\documentclass[12pt]{amsart}
\usepackage{amsmath, amssymb, latexsym, amsthm}
\usepackage{mathrsfs}
\usepackage[all]{xy}

\newcommand{\Bc}[9]{\bibitem{#1} {#2}, \emph{#3}, in: \textbf{#4} (#5), #6 #7, #8--#9.}

\newcommand{\Pa}[9]{\bibitem{#1} {#2}, \emph{#3}, {#4} \textbf{#5} ({#6}), {#7}--{#8}.}
\newcommand{\ed}{

\end{document}
}

      \newenvironment{changemargin}[2]{\begin{list}{}{
         \setlength{\topsep}{0pt}\setlength{\leftmargin}{0pt}
         \setlength{\rightmargin}{0pt}
         \setlength{\listparindent}{\parindent}
         \setlength{\itemindent}{\parindent}
         \setlength{\parsep}{0pt plus 1pt}
         \addtolength{\leftmargin}{#1}\addtolength{\rightmargin}{#2}
         }\item }{\end{list}}

\newcommand{\intvl}[2]{{[#1(#2),\allowbreak #1(#2\!+\!1))}}

\newcommand{\Dfin}{\mathfrak{D}_\mathrm{fin}}

\newcommand{\arx}[1]{\texttt{http://arxiv.org/abs/#1}}
\newcommand{\bq}{\begin{quote}}
\newcommand{\eq}{\end{quote}}

\newcommand{\CH}{the Continuum Hypothesis}

\newcommand{\Cantor}{{\{0,1\}^\N}}

\newcommand{\N}{\mathbb{N}}
\newcommand{\NN}{{\N^{\N}}}

\newcommand{\PN}{{P(\N)}}
\newcommand{\roth}{{[\N]^{\aleph_0}}}

\newcommand{\seq}[1]{\{#1\}_{n\in\N}}
\newcommand{\setseq}[1]{\{#1 : n\in\N\}}

\newcommand{\op}{\operatorname}

\newcommand{\maxfin}{\op{maxfin}}

\newcommand{\cI}{\mathcal{I}}

\newcommand{\scrA}{\mathscr{A}}
\newcommand{\scrB}{\mathscr{B}}

\newcommand{\CG}{C_\Gamma}
\newcommand{\CL}{C_\Lambda}
\newcommand{\CT}{C_\Tau}

\newcommand{\CO}{C_\Omega}
\newcommand{\Tau}{\mathrm{T}}

\newcommand{\cF}{\mathcal{F}}

\newcommand{\cM}{\mathcal{M}}

\newcommand{\cO}{\mathcal{O}}

\newcommand{\cU}{\mathcal{U}}
\newcommand{\Union}{\bigcup}
\newcommand{\cV}{\mathcal{V}}
\newcommand{\cW}{\mathcal{W}}

\newcommand{\Impl}{\Rightarrow}
\long\def\forget#1\forgotten{}

\newcommand{\fb}{\mathfrak{b}}
\newcommand{\fc}{\mathfrak{c}}
\newcommand{\fd}{\mathfrak{d}}
\newcommand{\fg}{\mathfrak{g}}

\newcommand{\fu}{\mathfrak{u}}

\newcommand{\x}{\times}
\newcommand{\Iff}{\Leftrightarrow}
\newcommand\comp{^{\text{\tt c}}}
\newcommand{\nin}{\not\in}

\newcommand{\sbst}{\subseteq}

\newcommand{\sm}{\setminus}
\newcommand{\as}{\subseteq^*}

\renewcommand{\>}{\rangle}

\newcommand{\cov}{\mathsf{cov}}

\newtheorem{thm}{Theorem}[section]
\newtheorem{prop}[thm]{Proposition}

\newtheorem{dict}[thm]{Dictionary}
\newtheorem{prob}[thm]{Problem}
\newtheorem{lem}[thm]{Lemma}
\newtheorem{cor}[thm]{Corollary}

\theoremstyle{definition}
\newtheorem{defn}[thm]{Definition}
\theoremstyle{remark}
\newtheorem{rem}[thm]{Remark}

\newcommand{\be}{\begin{enumerate}}
\newcommand{\ee}{\end{enumerate}}
\newcommand{\bi}{\begin{itemize}}
\newcommand{\itm}{\item}
\newcommand{\ei}{\end{itemize}}

\forget
\forgotten

\newcommand{\sone}{\mathsf{S}_1}
\newcommand{\sfin}{\mathsf{S}_\mathrm{fin}}
\newcommand{\ufin}{\mathsf{U}_\mathrm{fin}}
\newcommand{\Split}{\mathsf{Split}}

\author{Boaz Tsaban}
\thanks{Supported by the Koshland Center for Basic Research.}
\address[Boaz Tsaban]{Department of Mathematics, Bar-Ilan University, Ra\-mat-Gan 52900, Israel;
and
Department of Mathematics, Weizmann Institute of Science, Rehovot 76100, Israel.}
\email{tsaban@math.biu.ac.il}
\urladdr{http://www.cs.biu.ac.il/\~{}tsaban}

\author{Lyubomyr Zdomskyy}
\address[Lyubomyr Zdomskyy]{Department of Mechanics and Mathematics, Iv\-an Franko Lviv National University, Universytetska 1, Lviv 79000, Ukraine;
and
Department of Mathematics, Weizmann Institute of Science, Rehovot 76100, Israel.}
\curraddr{Kurt G\"odel Research Center for Mathematical Logic, W\"ahringer Str.\ 25, A-1090 Vienna, Austria.}
\email{lzdomsky@gmail.com}

\title[Combinatorial images and semifilter trichotomy]{Combinatorial images of sets of reals and semifilter trichotomy}

\begin{document}

\begin{abstract}
Using a dictionary translating a variety of classical and modern
covering properties into combinatorial properties of continuous
images, we get a simple way to understand the interrelations
between these properties in ZFC and in the realm of the trichotomy
axiom for upward closed families of sets of natural numbers. While
it is now known that the answer to the Hurewicz 1927 problem is
positive, it is shown here that semifilter trichotomy implies a
negative answer to a slightly stronger form of this problem.
\end{abstract}


\maketitle

\section{Introduction and basic facts}

Unless otherwise indicated, all spaces considered here
are assumed to be separable, zero-dimensional, and metrizable.
Consequently, we may assume that all open covers are countable \cite{split}.
Since every such space is homeomorphic to a set of real numbers, our
results can be thought of as dealing with sets of reals.

\subsection{Covering properties}
Fix a space $X$.
An open cover $\cU$ of $X$ is
\emph{large} if each member of $X$ is contained in infinitely
many members of $\cU$.
$\cU$ is an \emph{$\omega$-cover} if $X$ is not in $\cU$ and for
each finite $F\sbst X$, there is $U\in\cU$ such that $F\subseteq U$.
$\cU$ is a \emph{$\gamma$-cover} of $X$ if it is infinite and for each $x\in X$,
$x$ is a member of all but finitely many members of $\cU$.

Let $\cO$, $\Lambda$, $\Omega$, and $\Gamma$ denote the collections of all countable open
covers, large covers, $\omega$-covers, and $\gamma$-covers of $X$, respectively.
Let
$\scrA$ and $\scrB$ be any of these classes. We consider the following three
properties which $X$ may or may not have.
\begin{itemize}
\item[$\sone(\scrA,\scrB)$:]
For each sequence $\seq{\cU_n}$ of members of $\scrA$,
there exist members $U_n\in\cU_n$, $n\in\N$, such that $\setseq{U_n}\in\scrB$.
\item[$\sfin(\scrA,\scrB)$:]
For each sequence $\seq{\cU_n}$
of members of $\scrA$, there exist finite
subsets $\cF_n\sbst\cU_n$, $n\in\N$, such that $\Union_{n\in\N}\cF_n\in\scrB$.
\item[$\ufin(\scrA,\scrB)$:]
For each sequence $\seq{\cU_n}$ of members of $\scrA$
which do not contain a finite subcover,
there exist finite subsets $\cF_n\sbst\cU_n$, $n\in\N$,
such that $\setseq{\cup\cF_n}\in\scrB$.
\end{itemize}

It was shown by Scheepers \cite{coc1} and by Just, Miller, Scheepers, and Szeptycki
\cite{coc2} that each of these properties, when
$\scrA,\scrB$ range over $\cO,\Lambda,\Omega,\Gamma$,
is either void or equivalent to one in the following
diagram (where an arrow denotes implication).
For these properties, $\cO$ can be replaced anywhere by $\Lambda$ without changing the property.

{\scriptsize
\begin{changemargin}{-4.5cm}{-3cm}
\begin{center}
$\xymatrix@R=10pt{
&
&
& \ufin(\cO,\Gamma)\ar[r]
& \ufin(\cO,\Omega)\ar[rr]
& & \sfin(\cO,\cO)
\\
&
&
& \sfin(\Gamma,\Omega)\ar[ur]
\\
& \sone(\Gamma,\Gamma)\ar[r]\ar[uurr]
& \sone(\Gamma,\Omega)\ar[rr]\ar[ur]
& & \sone(\Gamma,\cO)\ar[uurr]
\\
&
&
& \sfin(\Omega,\Omega)\ar'[u][uu]
\\
& \sone(\Omega,\Gamma)\ar[r]\ar[uu]
& \sone(\Omega,\Omega)\ar[uu]\ar[rr]\ar[ur]
& & \sone(\cO,\cO)\ar[uu]
}$
\end{center}
\end{changemargin}
}

$\sfin(\cO,\cO)$, $\ufin(\cO,\Gamma)$, $\sone(\cO,\cO)$
are the the classical properties of
Menger, Hurewicz, and Rothberger ($C'')$, respectively.
$\sone(\Omega,\Gamma)$ is the Gerlits-Nagy $\gamma$-property.
Additional properties in the diagram were studied by Arkhangel'ski\v{i},
Sakai, and others. Some of the properties are relatively new.

We also consider the following type of properties.
\bi
\itm[$\Split(\scrA,\scrB)$:] Every cover $\cU\in\scrA$ can be split
into two disjoint subcovers $\cV$ and $\cW$ which contain elements of $\scrB$.
\ei
Here too, letting $\scrA,\scrB\in\{\Lambda,\Omega,\Gamma\}$
we get that some of the properties are trivial and several
equivalences hold among the remaining ones.
The surviving properties are
$$
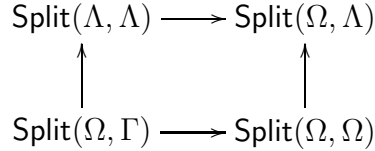
\xymatrix{
\Split(\Lambda, \Lambda) \ar[r] & \Split(\Omega, \Lambda)\\
\Split(\Omega, \Gamma) \ar[u]\ar[r] & \Split(\Omega, \Omega)\ar[u]\\
}$$
and no implication can be added to the diagram \cite{split}.
There are connections between the first and the second diagram,
e.g., $\Split(\Omega, \Gamma)=\sone(\Omega,\Gamma)$ \cite{split}, and
both $\ufin(\cO,\Gamma)$ and $\sone(\cO,\cO)$ imply $\Split(\Lambda,\Lambda)$.
Similarly, Scheepers proved that $\sone(\Omega,\Omega)$ implies $\Split(\Omega,\Omega)$ \cite{coc1}.

Let $C$, $\CL$, $\CO$, and $\CG$ denote the collections of all countable \emph{clopen}
covers, large covers, $\omega$-covers, and $\gamma$-covers of $X$, respectively.

It is often the case that we do not get anything new if we replace
an ordered pair of families of open covers
by the corresponding ordered pair of families of clopen covers.
However, some problems remain open.

\begin{prob}
Is any of the properties
\be
\itm $\sfin(\Gamma,\Omega)$, $\sone(\Gamma,\Gamma)$, $\sone(\Gamma,\Omega)$, $\sone(\Gamma,\cO)$;
\itm $\Split(\Lambda, \Lambda)$, $\Split(\Omega, \Lambda)$, $\Split(\Omega, \Omega)$;
\ee
equivalent to the corresponding property for clopen covers?
\end{prob}

In any case, the clopen version of each property is formally weaker.

\subsection{Combinatorial images}
The \emph{Baire space} $\NN$ and the \emph{Cantor space} $\Cantor$ are both equipped with the product topology.
$\PN$, the collection of all subsets of $\N$, is identified with
$\Cantor$ via characteristic functions, and inherits its topology.
The \emph{Rothberger space} $\roth$, consisting of all infinite sets of natural numbers,
is a subspace of $\PN$ and is homeomorphic to $\NN$.

For $a,b\in\roth$, $a$ is an \emph{almost subset} of $b$, $a\as b$, if
$a\sm b$ is finite.

\begin{defn}
A \emph{semifilter} is a nonempty family $F\sbst\roth$
containing all almost-supersets of its elements.
For a nonempty family $S\sbst\roth$,
$$\<S\>=\{b\in\roth : (\exists a\in S)\ a\as b\}$$
is the semifilter generated by $S$.
If $F=\<S\>$, then we say that $S$ is a \emph{base} for $F$.
A \emph{filter} is a semifilter closed under finite intersections,
and a \emph{subbase} for a filter is a family which, after closing under
finite intersections, becomes a base for that filter.
\end{defn}

The names of the combinatorial notions in the following dictionary are standard,
and a good reference for these is Blass' \cite{BlassHBK}.
We say that $g\in\NN$ is a \emph{guessing function} for $Y\sbst\NN$ if for each
$f\in Y$, $g(n)=f(n)$ for infinitely many $n$. In this case, we say that $Y$ is \emph{guessable}.
The following will be used throughout the paper without further notice.

\begin{dict}\label{dict}
The \emph{negation} of each property in the left column of the following table is equivalent to
having a continuous image in the relevant space ($\NN$ in the first block, and $\roth$ in the second) with
the corresponding property in the right column.
\begin{center}
\begin{tabular}{|l|lr|}
\hline
$\sfin(\cO,\cO)$    & dominating          & \cite{Rec94}\\
\hline
$\ufin(\cO,\Gamma)$ & unbounded           & \cite{Rec94}\\
\hline
$\sone(\cO,\cO)$    & not guessable       & \cite{Rec94}\\
\hline
$\ufin(\cO,\Omega)$ & finitely-dominating & \cite{huremen1}\\
\hline
\hline
$\Split(\CL,\CL)$ & reaping & \cite{split}\\
\hline
$\Split(\CO,\CL)$ & ultrafilter base & \cite{split}\\
\hline
$\Split(\CO,\CO)$ & ultrafilter subbase & \cite{split}\\
\hline
$\Split(\CT,\CT)$ & simple $P$-point base & \cite{split}\\
\hline
\end{tabular}
\end{center}
The analogous assertions for countable Borel covers, with ``continuous'' replaced
by ``Borel'', also hold \cite{CBC, split}.
\end{dict}

\subsection{Semifilter trichotomy, reformulated}
We now define one of the paper's main tools.
Recall that the \emph{Fr\'echet filter} is the set of all
cofinite subsets of $\N$.
\begin{defn}
For $a\in\roth$ and an increasing $h\in\NN$, define
$$a/h = \{n : a\cap \intvl{h}{n}\neq\emptyset\}.$$
For $S\sbst\roth$, define $S/h = \{a/h : a\in S\}$.
\emph{semifilter trichotomy} is the statement:
For each semifilter $S$, there is an increasing $h\in\NN$
such that $S/h$ is either the Fr\'echet filter, or an ultrafilter, or $\roth$.
\end{defn}

\begin{rem}
Semifilter trichotomy is consistent:
Blass and Laflamme \cite{BlassLaf89}, using a model invented for another purpose
in Blass and Shelah \cite{BlSh87}, proved
that the inequality $\fu<\fg$, where $\fu$ is the ultrafilter number and $\fg$ is the groupwise density number,
is consistent. Laflamme \cite{Laf92} proved that semifilter trichotomy follows from $\fu<\fg$.

In fact, Blass proved that semifilter trichotomy also implies $\fu<\fg$ \cite{Blass90}, and thus
semifilter trichotomy is equivalent to $\fu<\fg$.
\end{rem}

When speaking of an element $a\in\roth$ as an element of $\NN$, we do this by
identifying $a$ with its increasing enumeration. This identification gives
a homeomorphism from $\roth$ onto the set of increasing elements in $\NN$.
Thus, we say that a family $S\sbst\roth$ is \emph{unbounded} if it is unbounded
when viewed as a subset of $\NN$.

\begin{defn}
An increasing $h\in\NN$ is a \emph{(flat) slalom} for a family $S\sbst\roth$ if
for each $a\in S$, for all but finitely many $n$, $a\cap\intvl{h}{n}\neq\emptyset$.
\end{defn}

It is easy to see (e.g., \cite{hureslaloms}) that $S$ has a slalom if, and only if, it is bounded.

\begin{cor}\label{slalom}
A family $S\sbst\roth$ is bounded if, and only if, there is an increasing $h\in\NN$ such that $\<S/h\>$
is the Fr\'echet filter.
\end{cor}
\begin{proof}
$\<S/h\>$ is the Fr\'echet filter if, and only if, for each $a\in S$, $a/h$ is cofinite,
that is, $h$ a slalom for $S$.
\end{proof}

\begin{thm}\label{tricho}
The following assertions are equivalent:
\be
\itm Semifilter trichotomy.
\itm For each unbounded $S\sbst\roth$, there is an increasing $h\in\NN$ such that
$S/h$ is a base for either an ultrafilter, or for $\roth$.
\itm For each unbounded $S\sbst\roth$, there is an increasing $h\in\NN$ such that $S/h$ is reaping.
\ee
\end{thm}
\begin{proof}
$(1\Iff 2)$ $S/h$ is always a base for $\<S\>/h$. Use Corollary \ref{slalom}.

$(2\Impl 3)$ Is trivial.

$(3\Impl 1)$ Each intersection of two unbounded semifilters is unbounded \cite{BlassHBK}.
Let $S$ be a semifilter, and assume that for each $h$, $S/h\neq\roth$ and is
not the Fr\'echet filter.
Then the same is true for $S^+=\{a\in\roth : a\comp\nin S\}$.
Let $U$ be an ultrafilter. As $S^+,U$ are unbounded, $F=S^+\cap U$ is unbounded.
Thus, there is $h$ such that the semifilter $F/h$ is reaping.
As $F/h$ is a reaping subset of an ultrafilter $U/h$, $F/h=U/h$.
It follows that $U/h\sbst S^+/h$, and as $U/h$ is an ultrafilter, we have that
$S/h = (S^+/h)^+\sbst (U/h)^+ = U/h$ is an ultrafilter.
\end{proof}

\section{Warm up: Three basic results in ZFC}

The results below were originally proved using sophisticated manipulations of open covers.
The combinatorial proofs given here are direct generalizations of arguments
from the theory of cardinal characteristics of the continuum.

\begin{thm}[Scheepers \cite{coc1}]\label{split1}
$\ufin(\cO,\Gamma)$ implies $\Split(\CL,\CL)$.
\end{thm}
\begin{proof}
Assume that $Y\sbst\roth$ is a continuous image of $X$.
As $X$ has the Hurewicz property, $Y$ has a slalom $h$ \cite{hureslaloms}.
It suffices to show that $Y$ is not reaping.
Indeed, let $a=\Union_n\intvl{h}{2n}$.
Then for each $y\in Y$, both $y\cap a$ and $y\cap a\comp$ are infinite.
\end{proof}

\begin{thm}[Scheepers \cite{coc1}]\label{rothsplit}
$\sone(\cO,\cO)$ implies $\Split(\CL,\CL)$.
\end{thm}
\begin{proof}
Assume that $X$ satisfies $\sone(\cO,\cO)$,
and $Y\sbst\roth$ is a continuous image of $X$.
For each $y\in Y$, define $f_y\in \prod_n[\N]^{2n}$
by $f_y(n) = \{y(1),\dots,y(2n)\}$.

For each $n$, we can identify $[\N]^{2n}$ with $\N$ and therefore
identify $\prod_n[\N]^{2n}$ with $\NN$ in a natural way.
$Z=\{f_y : y\in Y\}$ is a continuous image of $Y$, and thus
there is a guessing function $g\in\prod_n[\N]^{2n}$
for $Z$.
For each $n$, let $i_n,j_n$ be distinct members of $g(n)\sm \{i_1,\dots,i_{n-1},j_1,\dots,j_{n-1}\}$.
Take $I=\setseq{i_n}, J=\setseq{j_n}$.

For each $y\in Y$ there are infinitely many $n$ such that $g(n)=f_y(n)$,
and therefore both $I\cap y$ and $J\cap y$ are infinite.
As $I\cap J=\emptyset$, $Y$ is not reaping.
\end{proof}

Scheepers proved in \cite{coc1} that $\sone(\Omega,\Omega)$ implies $\Split(\Omega,\Omega)$.
Ko\v{c}inac and Scheepers \cite{coc7} proved that if all finite powers of $X$ satisfy $\ufin(\cO,\Gamma)$, then
$X$ satisfies $\Split(\Omega,\Omega)$. Both results are generalized in a single result
from \cite{split}, asserting that if all finite powers of $X$ satisfy $\Split(\Omega,\Lambda)$,
then $X$ satisfies $\Split(\Omega,\Omega)$. The same proof works in the clopen case, but it is
quite complicated. We give a simple proof.

\begin{thm}[\cite{split}]
If all finite powers of $X$ satisfy $\Split(\CO,\CL)$, then
$X$ satisfies $\Split(\CO,\CO)$.
\end{thm}
\begin{proof}
Assume that $X$ does not satisfy $\Split(\CO,\CO)$,
and let $Y\sbst\roth$ be a continuous image of $X$ which is
a subbase for an ultrafilter.
Note that all finite powers of $Y$ satisfy $\Split(\CO,\CL)$.
For each $k$, define
$\Psi_k:Y^k\to\roth$ by
$$(a_1,\dots,a_k)\mapsto a_1\cap\dots\cap a_k$$
for each $a_1,\dots,a_k\in Y$.
$\Psi_k$ is continuous, and therefore its image satisfies $\Split(\CO,\CL)$.
As $\Split(\CO,\CL)$ is $\sigma$-additive \cite{split},
$Z=\Union_k\Psi_k[Y^k]$ satisfies $\Split(\CO,\CL)$,
and $Z$ is a base for an ultrafilter -- a contradiction.
\end{proof}

\section{When semifilter trichotomy holds}

The second part of the following theorem was
proved in \cite{SF1}, using much more complicated arguments.

\begin{thm}\label{split2}
Assume semifilter trichotomy.
Then
$$\ufin(\cO,\Gamma)=\Split(\CL,\CL).$$
In particular, $\ufin(\cO,\Gamma)=\Split(\Lambda,\Lambda)$.
\end{thm}
\begin{proof}
By Theorem \ref{split1}, it suffices to prove that every space $X$
satisfying $\Split(\CL,\CL)$, satisfies $\ufin(\cO,\Gamma)$.

Indeed, assume that a continuous image $Y\sbst\roth$ of $X$ is unbounded.
By Lemma \ref{tricho}, there is an increasing $h\in\NN$ such that
$Y/h$ (a continuous image of $Y$, and therefore of $X$) is
reaping. Thus, $X$ does not satisfy $\Split(\CL,\CL)$.

For the last assertion of the theorem, use Scheepers' result
that $\ufin(\cO,\Gamma)$ implies $\Split(\Lambda,\Lambda)$ \cite{coc1},
and the trivial fact that $\Split(\Lambda,\Lambda)$ implies $\Split(\CL,\CL)$.
\end{proof}

The following natural concept, due to Ko\v{c}inac and Scheepers \cite{coc7},
will appear several times in this paper. We introduce it using the
self-explanatory terminology of \cite{GlCovs}.

\begin{defn}
A cover $\cU$ of $X$ is \emph{$\gamma$-glueable} if
$\cU$ can be partitioned into infinitely many finite pieces,
such that either each piece covers $X$, or else the unions of
the pieces form a $\gamma$-cover of $X$.
$\gimel(\Gamma)$ is the family of all open $\gamma$-glueable covers of $X$.
\end{defn}

The Gerlits-Nagy property $(*)$ is defined in \cite{GN}. In
\cite{coc7} it is shown that this property is equivalent to
$\sone(\Lambda,\gimel(\Gamma))$.

\begin{cor}
Assume semifilter trichotomy. Then
$$\sone(\Lambda,\gimel(\Gamma))=\sone(\cO,\allowbreak\cO).$$
\end{cor}
\begin{proof}
$\sone(\Lambda,\gimel(\Gamma))= \ufin(\cO,\Gamma)\cap\sone(\cO,\cO)$ \cite{coc7}.
Apply Theorems \ref{rothsplit} and \ref{split2}.
\end{proof}

A classical problem of Hurewicz asks whether $\ufin(\cO,\Gamma)\neq \sfin(\cO,\cO)$.
Chaber and Pol \cite{ChaPol} gave a positive answer outright in ZFC (see \cite{SFH}).
However, we can show that a slightly stronger assertion is consistently true.
The property $\Split(\Omega,\Lambda)$ is not very restrictive:
E.g., it holds for every analytic space \cite{split}.

\begin{thm}\label{split3}
Assume semifilter trichotomy.
Then
$$\ufin(\cO,\Gamma)=\sfin(\cO,\cO)\cap\Split(\CO,\CL).$$
In particular, $\ufin(\cO,\Gamma)=\sfin(\cO,\cO)\cap\Split(\Omega,\Lambda)$.\hfill\qed
\end{thm}
\begin{proof}
Any base for $\roth$, when viewed as a subset of $\NN$, is dominating.
Thus, the proof is the same as in Theorem \ref{split2}.
\end{proof}

\begin{rem}
Theorem \ref{split3} cannot be improved to get $\ufin(\cO,\Gamma)=\Split(\Omega,\allowbreak\Lambda)$
from semifilter trichotomy,
since any analytic set (in particular, $\NN$) satisfies $\Split(\Omega,\Lambda)$ \cite{split}.
Moreover, some axiom is necessary to get the equality in Theorem \ref{split3},
since even the stronger property $\sone(\Omega,\Omega)$ does not imply $\ufin(\cO,\Gamma)$ \cite{coc2}.
\end{rem}

\begin{rem}
In \cite{SF1}, a space $X$ is called \emph{almost Menger} if
for each large open cover $\setseq{U_n}$ of $X$, setting
$Y = \{\{n : x\in U_n\} : x\in X\}$ we have that
for each increasing $h\in\NN$, $Y/h$ is not a base for $\roth$.
It is shown there that if $X$ satisfies $\sfin(\cO,\cO)$ then $X$ is
almost Menger, and we are asked whether the converse holds.
As a base for $\roth$ must have cardinality $\fc$, we have that
the answer is negative when $\fd<\fc$.

On the other hand, the proof of Theorem \ref{split3} shows that assuming semifilter trichotomy,
if $X$ is almost Menger
and satisfies $\Split(\Omega,\Lambda)$, then $X$ satisfies $\ufin(\cO,\allowbreak\Gamma)$.
\end{rem}

We now give a simple proof for the following result, which involves no splitting properties.

\begin{thm}[\cite{SF1}]\label{split4}
Assume semifilter trichotomy. Then
$$\ufin(\cO,\Omega)=\sfin(\cO,\cO).$$
\end{thm}
\begin{proof}
Assume that $X$ satisfies $\sfin(\cO,\cO)$, and that $Y\sbst\NN$ is a continuous image of $X$.
We may assume that all elements in $Y$ are increasing.
$Y$ is not dominating. Choose an increasing $g\in\NN$ witnessing that.
The collection $Z$ of the sets $[f\le g]=\{n : f(n)\le g(n)\}$, $f\in Y$, is a continuous image
of $Y$ in $\roth$. Thus, for each increasing $h\in\NN$, $Z/h$ is not a base for $\roth$.
By semifilter trichotomy, there is an increasing $h\in\NN$ such that
$Z/h$ is a base for a filter $F$ ($F$ is either an ultrafilter or the Fr\'echet filter).
We will show that $Y$ is bounded with respect to $F$.

Indeed, define $\tilde g\in\NN$ by $\tilde g(n)=g(h(n+1))$ for all $n$.
For each $f\in Y$, let $a = [f\le g]/h\in F$.
For each $n\in a$, choose $k\in [f\le g]\cap\intvl{h}{n}$.
Then
$$f(n)\le f(h(n))\le f(k)\le g(k)\le g(h(n+1))=\tilde g(n).$$
Thus, $a\sbst [f\le \tilde g]$. As $a\in F$, $[f\le \tilde g]\in F$.
As $F$ is a filter, $\tilde g$ witnesses that $Y$ is not finitely dominating.
\end{proof}

We have thus obtained a simple proof for the following.

\begin{cor}[\cite{huremen2}]
Assume semifilter trichotomy. Then $\ufin(\cO,\Omega)$ is $\sigma$-additive.\hfill\qed
\end{cor}

\section{$\ufin(\cO,\Omega)$ revisited}

Now that we know that consistently $\ufin(\cO,\Omega)=\sfin(\cO,\cO)$,
we can step back to ZFC and ask whether some nontrivial properties of
$\sfin(\cO,\cO)$ can be transferred to $\ufin(\cO,\Omega)$.
This is the purpose of this section.

In \cite{MGD} it is proved that if $X$ satisfies $\sfin(\cO,\cO)$,
then for each continuous image $Y$ of $X$ in $\NN$, the set
$$G=\{g\in\NN : (\forall f\in Y)\ g\not\le^* f\}$$
is nonmeager. In particular, this is true for
$\ufin(\cO,\Omega)$, but this is not the correct assertion for
that property.
For $Y\sbst\NN$, let
$$\maxfin(Y)=\{\max\{f_1,\dots,f_k\} : k\in\N,\ f_1,\dots,f_k\in Y\}.$$
Then $X$ satisfies $\ufin(\cO,\Omega)$ if, and only if, for each continuous image $Y$ of $X$ in $\NN$,
$\maxfin(Y)$ is not dominating.

\begin{thm}\label{manybounds}
For each space $X$, the following are equivalent.
\be
\itm $X$ satisfies $\ufin(\cO,\Omega)$.
\itm For each continuous image $Y$ of $X$ in $\NN$,
the set
$$G=\{g\in\NN : (\forall f\in\maxfin(Y))\ g\not\le^* f\}$$
is nonmeager.
\ee
\end{thm}
\begin{proof}
$(2\Impl 1)$ nonmeager sets are nonempty.

$(1\Impl 2)$ Assume that $X$ satisfies $\ufin(\cO,\Omega)$ and $Y\sbst\NN$ is a continuous image of $X$.
If $Y$ is bounded, then (2) holds trivially.
Assume that $Y$ is unbounded.
Let $g$ be a witness for the fact that $Y$ is not finitely dominating.
Take
$$Z=\{[f<g] : f\in Y\}.$$
$Z$ is a subbase for a filter. Extend this filter to a nonprincipal ultrafilter $F$.
For each $f\in Y$, $f\le_F g$. As $F$ is a filter, $\le_F$ is transitive,
so it suffices to show that the set
$$G'=\{f\in\NN : g\le_F f\}$$
is nonmeager. Since $F$ is a nonmeager semifilter, this is true \cite{SFH}.
(For an alternative approach see \cite{MGD} and Lemma 2.4 of Mildenberger, Shelah, and Tsaban \cite{MShT:847}.)
\end{proof}

The proof of Theorem \ref{manybounds} turned out easier than the corresponding one
for $\sfin(\cO,\cO)$. However, for $\sfin(\cO,\cO)$ we get slightly more:
If $X$ satisfies $\sfin(\cO,\cO)$, then for each continuous image $Y$ of $X$ in $\NN$,
the set
$$G=\{g\in\NN : (\exists f\in Y)\ g\le^* f\}$$
satisfies $\sfin(\cO,\cO)$ \cite{MGD}.
To see why this is indeed more, consider the following.

\begin{lem}
Assume that $Y$ is a subset of $\NN$ and satisfies
$\sfin(\cO,\cO)$. Then $Y$ is not comeager.
\end{lem}
\begin{proof}
Assume that $Y$ is comeager.
To each $f\in\NN$, assign the set
$$a_f=\setseq{f(0)+\dots+f(n)+n}.$$
$f\mapsto a_f$ is a homeomorphism from $\NN$ to $\roth$.
Thus, $Z=\{a_f : f\in Y\}$ satisfies $\sfin(\cO,\cO)$ and is comeager.
By a classical result of Talagrand \cite{CSF},
for each comeager subset $Z$ of $\roth$ there is
an increasing $h\in\NN$ such that $\<Z/h\>=\roth$.
It follows that $Z/h$ is dominating -- a contradiction.
\end{proof}

The following remains open.

\begin{prob}\label{Sch}
Assume that $X$ satisfies $\ufin(\cO,\Omega)$, and $Y\sbst\NN$ is a continuous image of $X$.
Does it follow that
$$G=\{g\in\NN : (\exists k)(\exists f_1,\dots,f_k\in Y)\ g\le^* \max\{f_1,\dots,f_k\}\}$$
satisfies $\ufin(\cO,\Omega)$?
\end{prob}

In the remainder of this section we will show that the auxiliary results proved
in \cite{MGD} for $\sfin(\cO,\cO)$, which are interesting in their own right,
also hold for $\ufin(\cO,\Omega)$.

\medskip

It is consistent that $\ufin(\cO,\Omega)$ is not even preserved under taking finite unions.
In fact, this follows from \CH{} (or even just $\cov(\cM)=\fc$) \cite{huremen2}.
However, something is still provable about unions of spaces satisfying $\ufin(\cO,\Omega)$.
Let $\cov(\Dfin)$ denote the minimal cardinality of a partition of $\NN$ into families
which are not finitely dominating. This is the same as the
minimal cardinality of a partition of any dominating family in $\NN$ into families
which are not finitely dominating. $\max\{\fb,\fg\}\le\cov(\Dfin)$, and it is consistent
that strict inequality holds \cite{MShT:847}.

\begin{prop}\label{unions}
Assume that $Z$ is a space, and $\cI\sbst P(Z)$ satisfies:
\be
\itm For each finite $\cF\sbst\cI$, there is $X\in\cI$ such that $\cup\cF\sbst X$;
\itm Each $X\in\cI$ satisfies $\ufin(\cO,\Omega)$;
\itm $|\cI|< \cov(\Dfin)$.
\ee
Then $\cup\cI$ satisfies $\ufin(\cO,\Omega)$.
\end{prop}
\begin{proof}
Assume that $\Psi:\cup\cI\to \NN$ is continuous.
By (2), for each $X\in\cI$, $\Psi[X]$ is not finitely dominating,
and therefore $\maxfin(\Psi[X])$ is not finitely dominating.
By (1),
$$\maxfin(\Psi[\cup\cI])=\Union_{X\in\cI}\maxfin(\Psi[X]).$$
By (3), $\maxfin(\Psi[\cup\cI])$ is not dominating, that is,
$\Psi[\cup\cI]$ is not finitely dominating.
\end{proof}

As $\ufin(\cO,\Omega)$ is hereditary for closed subsets,
Proposition \ref{unions} implies the following.

\begin{cor}\label{FsHered}
$\ufin(\cO,\Omega)$ is hereditary for $F_\sigma$ subsets.\hfill\qed
\end{cor}

Another interesting corollary is the following.

\begin{cor}\label{inc}
$\ufin(\cO,\Omega)$ is preserved under taking countable increasing unions.\hfill\qed
\end{cor}

Finally, we have the following.

\begin{prop}\label{Ksigma}
Assume that $X$ satisfies $\ufin(\cO,\Omega)$ and $K$ is $\sigma$-compact.
Then $X\x K$ satisfies $\ufin(\cO,\Omega)$.
\end{prop}
\begin{proof}
By Corollary \ref{inc}, we may assume that $K$ is compact
(one can also manage without that).
Assume that $\cU_1,\cU_2,\dots$, are countable open covers of $X\x K$.
For each $n$, enumerate $\cU_n=\{U^n_m : m\in\N\}$.
For each $n$ and $m$ set
$$V^n_m = \left\{x\in X : \{x\}\x K\sbst \Union_{k\le m}U^n_k\right\}.$$
Then $\cV_n = \{V^n_m : m\in\N\}$ is an open cover of $X$.
As $X$ satisfies $\ufin(\cO,\Omega)$, we can choose for each $n$ an $m_n$
such that for each finite $F\sbst X$, there is $n$ such that $F\sbst \Union_{k\le m_n}V^n_k$.

Assume that $F\sbst X\x K$ is finite. Take finite $A\sbst X, B\sbst K$ such that
$F\sbst A\x B$. Let $n$ be such that $A\sbst \Union_{k\le m_n}V^n_k$.
Then for each $a\in A$, $a\x K\sbst \Union_{k\le m_n}U^n_k$, and therefore
$$A\x B\sbst A\x K\sbst \Union_{k\le m_n}U^n_k.\qedhere$$
\end{proof}

\begin{rem}
All properties in the Scheepers diagram are hereditary for \emph{closed} subsets.
As $\ufin(\cO,\Gamma)$, $\sfin(\cO,\cO)$, $\sone(\Gamma,\Gamma)$,
$\sone(\Gamma,\cO)$, and $\sone(\cO,\cO)$ are all $\sigma$-additive \cite{AddQuad},
they are all hereditary for $F_\sigma$ subsets.
Galvin and Miller \cite{GM} proved that
$\sone(\Omega,\Gamma)$ is also hereditary for $F_\sigma$ subsets.
$\sfin(\Omega,\Omega)$ is equivalent to satisfying $\sfin(\cO,\cO)$ in all finite
powers. As finite powers of $F_\sigma$ sets are $F_\sigma$,
$\sfin(\Omega,\Omega)$ is also hereditary for $F_\sigma$ subsets.
Similarly, $\sone(\Omega,\Omega)$ is equivalent to satisfying $\sone(\cO,\cO)$ in all finite
powers and is therefore also hereditary for $F_\sigma$ subsets.
By Corollary \ref{FsHered}, so is $\ufin(\cO,\Omega)$.
\end{rem}

\begin{prob}
Are $\sfin(\Gamma,\Omega)$ and $\sone(\Gamma,\Omega)$ hereditary for $F_\sigma$ subsets?
\end{prob}

\section{The revised Hurewicz Problem for general spaces}

As mentioned before, Theorem \ref{split3}
may be considered a consistent positive solution to a revised
version of the original Hurewicz Problem
(which had a negative solution in ZFC).

Since this result is new, we prove that it holds in general,
i.e., without any assumption on the spaces.

\begin{thm}
Assume semifilter trichotomy.
Then
$$\ufin(\cO,\Gamma)=\sfin(\cO,\cO)\cap\Split(\Omega,\Lambda)$$
for arbitrary topological spaces.
\end{thm}
\begin{proof}
Assume that $X$ satisfies $\sfin\allowbreak(\cO,\cO)\cap\Split(\Omega,\Lambda)$.
By $\sfin(\cO,\cO)$, we have that $X$ is Lindel\"of.
In \cite{coc7} it is proved that $\ufin(\cO,\Gamma)=\sfin(\Lambda,\gimel(\Gamma))$.\footnote{%
The proof in \cite{coc7} only requires that $X$ is Lindel\"of.}
As $\sfin(\cO,\cO)=\sfin(\Lambda,\Lambda)$ \cite{coc1, coc2}, we have
that for Lindel\"of spaces,
$$\ufin(\cO,\Gamma)=\sfin(\Lambda,\gimel(\Gamma))=\sfin(\Lambda,\Lambda)\cap\binom{\Lambda}{\gimel(\Gamma)}
= \sfin(\cO,\cO)\cap\binom{\Lambda}{\gimel(\Gamma)},$$
where $\binom{\Lambda}{\gimel(\Gamma)}$ means that every element of
$\Lambda$ contains an element of $\gimel(\Gamma)$.
It therefore remains to prove this latter property.

Let $\cU$ be a large open cover of $X$. As $X$ satisfies $\sfin(\Lambda,\Lambda)$, we may assume
that $\cU$ is countable and fix a bijective enumeration $\cU=\setseq{U_n}$.
Let
$$Y = \{\{n : x\in U_n\} : x\in X\}.$$
Choose an increasing $h\in\NN$ witnessing semifilter trichotomy for $\<Y\>$.
For each $n$, define
$$V_n=\Union_{k\in\intvl{h}{n}} U_k.$$

\subsubsection*{Case 1} There are infinitely many $n$ such that $V_n=X$.
Let $a\in\roth$ be the set of all these $n$.
Taking $g(0)=0$ and $g(n)=h(a(n-1))$ for $n>0$, we have that
the sets $\cF_n=\{U_k : k\in\intvl{g}{n}\}$, $n\in\N$, form a partition
of $\cU$ showing that it is $\gamma$-glueable.

\subsubsection*{Case 2} There are only finitely many $n$ such that $V_n=X$.
Removing finitely many elements from $\cU$, we may assume that there are no such $n$.
(We can add these elements later to one of the pieces of the partition).

Assume that $Y/h$ is a base for an ultrafilter.
Then for each finite $a_1,\dots,a_k\in X$, there is $n\in a_1/h\cap\dots a_k/h$,
that is, $a_1,\dots,a_k\in V_n$. Thus, $\cV=\setseq{V_n}$
is an open $\omega$-cover of $X$.
As $Y/h$ is reaping, $\cV$ cannot be split into two large covers of $X$.
This contradicts $\Split(\Omega,\Lambda)$.

As $Y$ satisfies $\sfin(\cO,\cO)$, $Y/h$ is not a base for $\roth$ \cite{SF1}.

If follows that all elements in $Y/h$ are cofinite, that is,
for each $x\in X$ and all but finitely many $n$, $x\in V_n$.
This shows that $\cU$ is $\gamma$-glueable.
\end{proof}

It is not always the case that theorems of the discussed sort can be transferred
from sets of reals to arbitrary spaces.
We conclude the paper with an example for that.

It is known that for sets of reals, $\ufin(\cO,\Gamma)=\binom{\Lambda}{\gimel(\Gamma)}$ \cite{hureslaloms}.
Had we been able to prove this for general topological spaces, this would
have made the last proof shorter. Unfortunately, this can be refuted in a strong sense.

\begin{prop}
There exists a hereditarily Lindel\"of $T_1$ space $S$
satisfying $\binom{\Lambda}{\gimel(\Gamma)}$, but
not even $\sfin(\cO,\cO)$.
\end{prop}
\begin{proof}
Consider the topology $\tau$ on $\N$ generated by the sets $\setseq{[0,n)}$.
$\tau$ gives a product topology $\nu$ on $\NN$.
$(\NN,\nu)$ does not satisfy $\sfin(\cO,\cO)$: Indeed, consider the
open covers $\cU_n = \{U^n_m : m\in\N\}$ with $U^n_m = \{f\in\NN : f(n)\le m\}$.

Let $\mu$ be the topology generated by $\{U\setminus A: U\in\nu, A\sbst\NN\mbox{ is finite}\}$
as a base, and take $S=(\NN,\mu)$.
Clearly, $S$ is $T_1$.
As $\nu\sbst\mu$, $S$ does not satisfy $\sfin(\cO,\cO)$.
As $\mu$ is contained in the standard product topology on $\NN$,
$S$ is hereditarily Lindel\"of.

Assume that $\cU\sbst\mu$ is a large cover of $\NN$.
As $(\NN,\mu)$ is hereditarily Lindel\"of, we may assume that
$\cU$ is countable \cite{split}, and enumerate it bijectively as $\cU=\setseq{U_n\sm F_n}$,
where each $U_n\in\nu$ and each $F_n$ is a finite subset of $\NN$.
Let $D=\Union_n F_n$.
For a sequence $\cF=\setseq{X_n}$, and $f\in\NN$, write $f_\cF = \{n : f\in X_n\}$.

For each finite $F\sbst\NN$ let $g=\max F$.
Let $n$ be such that $g\in U_n\sm F_n$. Then $F\sbst U_n$.
It follows that $\cV=\setseq{U_n}$ is an $\omega$-cover of $\NN$
by sets open in the standard topology on $\NN$.
Consequently, $\cV$ is a $\gamma$-glueable cover of $\NN$ (Sakai \cite{Sakai06a}).
Then $\{f_\cV : f\in\NN\sm D\}$ is bounded.
Note that for each $f\nin D$, $f_\cV= f_\cU$, and therefore
$\{f_\cU : f\in\NN\sm D\}$ is bounded.
As $D$ is countable, $\{f_\cU : f\in D\}$ is also bounded,
and therefore $\{f_\cU : f\in \NN\}$ is bounded, that is,
$\cU$ is $\gamma$-glueable.
\end{proof}

\ed